\documentclass[oneside,10pt]{article}          
\usepackage[b5paper]{geometry}	    
\usepackage{amsfonts,amsmath,latexsym,amssymb} 
\usepackage{theorem}                
\usepackage{mathrsfs,upref}         
\usepackage{mathptmx}		    
\usepackage{elemath2024}	            
\usepackage{xcolor}
\usepackage[super]{nth}

\newtheorem{thm}{Theorem}[section]

\newtheorem{cor}[thm]{Corollary}

\theoremstyle{definition}
\newtheorem{rem}[thm]{Remark}

\newtheorem{defi}[thm]{Definition}

\numberwithin{equation}{section}       

\def\esup{\operatornamewithlimits{ess\,sup}}

\def\Id{\operatorname{\rm I}}

\def\Cop{\operatorname{Cop}}

\def\esup{\operatornamewithlimits{ess\,sup}}
\def\Id{\operatorname{I}}

\def\esup{\operatornamewithlimits{ess\,sup}}

\def\Id{\operatorname{I}}

\def\Ces{\operatorname{Ces}}
\def\Cop{\operatorname{Cop}}

\def\qq{\qquad}

\def\hra{\hookrightarrow}

\def\la{\lambda}

\def\I{(0,\infty)}

\def\R{\mathbb R}

\def\R{\mathbb R}

\def\L{{\mathcal L}}
\def\W{{\mathcal W}}

\def\mp{{\mathfrak M}}

\def\O{\Omega}

\def\la{\lambda}

\def\I{(a,b)}

\allowdisplaybreaks

\begin{document}
	
	\title[On weighted Ces\`{a}ro function spaces]{On weighted Ces\`{a}ro function spaces}
	
	
	\author{Amiran Gogatishvili and   Tu\u{g}\c{c}e \"{U}nver}
	
	\address{A. Gogatishvili,
		Institute of Mathematical of the  Czech Academy of Sciences, \v Zitn\'a~25, 115~67 Praha~1, Czech Republic\\
		\email{gogatish@math.cas.cz}}
	
	\address{T. \"{U}nver, Faculty of Engineering and Natural Sciences, Department of Mathematics,
		Kirikkale University, 71450, Yahsihan, Kirikkale, T\"{u}rkiye \\
		\email{tugceunver@kku.edu.tr}}

	\CorrespondingAuthor{T. \"{U}nver}
	
	\dedicated{Dedicated to  Professor Lars-Erik Persson on the occasion of the \nth{80} anniversary of his birthday}      
	
	\date{20.09.2024}                               
	
	
	
	\thanks{The research of Amiran Gogatishvili was partially supported by the Czech Academy of Sciences (RVO 67985840), by the Czech Science Foundation (GA\v{C}R), grant no: 23-04720S, and by  Shota Rustaveli National Science Foundation of Georgia (SRNSFG), grant no: FR22-17770. The research of Tu\u{g}\c{c}e \"{U}nver was supported by the Grant of The Scientific and Technological Research Council of Turkey (TUBITAK), Grant No: 1059B192000075.
	}

	\begin{abstract}
		The main objective of this paper is to provide a comprehensive demonstration of recent results regarding the structures of the weighted Cesàro and Copson function spaces. These spaces' definitions involve local and global weighted Lebesgue norms; in other words, the norms of these spaces are generated by positive sublinear operators and by weighted Lebesgue norms. The weighted Lebesgue spaces are the special cases of these spaces with a specific set of parameters.  
		
		Our primary method of investigating these spaces will be the so-called discretization technique. Our technique will be the development of the approach initiated by K.G. Grosse-Erdmann, which allows us to obtain the characterization in previously unavailable situations, thereby addressing decades-old open problems. 
		
		We investigate the relation  (embeddings) between weighted Cesàro and Copson function spaces. The characterization of these embeddings can be used to tackle the problems of characterizing pointwise multipliers between weighted Cesàro and Copson function spaces, the characterizations of the associate spaces of Cesàro (Copson) function spaces, as well as the relations between local Morrey-type spaces. 
	\end{abstract}
	
	\maketitle

	

	\section{Introduction and basic definitions}\label{introduction}


	Throughout the paper, $0\leq a \leq b \leq \infty$. By ${\mathcal M} (I)$ we denote the set of all
	measurable functions on $I$. The symbol ${\mathcal M}^+ (I)$ stands
	for the collection of all $f\in{\mathcal M} (I)$ which are
	non-negative on $I$. The family of all weight functions (also called just weights) on $I$, that is, measurable, positive, and finite, a.e. on $I$, is given by $\W (I)$.  
	
	We always denote $c$ and $C$ as positive constants, which are independent of main parameters but may differ from line to line. However, a constant with subscript or superscript such as $c_1$ does not change in different occurrences. By $A\lesssim B$ ($B\gtrsim A$) we mean that $A\leq \la B$, where $\la>0$ depends on inessential parameters. If $A\lesssim B$ and $B\lesssim A$, we write $A\approx B$ and say that $A$ and $B$ are equivalent.  We will denote by $\bf 1$ the function ${\bf 1}(x) = 1$, $x \in \R$.  
	Since the expressions on our main results are too
	long, to make the formulas plain, we sometimes omit the differential
	element $dx$. 
	
	We put $1/ (\pm \infty) =0$, $0/0 =0$, $0\cdot (\pm \infty)=0$. Moreover, we omit arguments of integrands as well as differentials in integrals when appropriate in order to keep the expository as short as possible.
	
	For $p\in (0,\infty]$ and $w\in \W(I)$, we define the functional $\|\cdot\|_{p,w,I}$ on $\mathcal M(I)$ by 
	\begin{equation*}
		\|f\|_{p,w,I} : = \left\{
		\begin{array}{cl}
			\bigg(\int_I |f(x)|^p w(x)^p\,dx \bigg)^{\frac{1}{p}} & \qq\mbox{if}\qq p<\infty, \\
			\esup\limits_{x\in I} |f(x)|w(x) & \qq\mbox{if}\qq p=\infty.
		\end{array}
		\right.
	\end{equation*}
	The weighted Lebesgue space $L_p(w,I)$ is given by
	\begin{equation*}
		L_p(w,I) \equiv L_{p,w}(I) : = \{f\in \mp (I):\,\, \|f\|_{p,w,I}  < \infty\},
	\end{equation*}
	and it is equipped with the quasi-norm $\|\cdot\|_{p,w,I}$. When $I = \I$, we often write simply $L_{p,w}$ and $L_p(w)$ instead of $L_{p,w}(I)$ and $L_p(w,I)$, respectively.
	
	The weighted Ces\`{a}ro and Copson function spaces are defined as
	follows:
	\begin{defi}\label{defi.2.1}
		Let  $0 <p, q \le \infty$, $u \in \mathcal M^+(a,b)$, $v\in \W(a,b)$. The weighted Ces\`{a}ro and Copson spaces are defined by 
		\begin{align*}
			\Ces_{p,q} (u,v) : & = \bigg\{ f \in \mathcal M^+(a,b): \|f\|_{\Ces_{p,q}(u,v)} : = \big\| \|f\|_{p,v,(a,\cdot)} \big\|_{q,u,(a,b)} < \infty \bigg\}, \\
			\intertext{and} 
			\Cop_{p,q}(u,v) : & = \bigg\{ f \in \mathcal M^+(a,b): \|f\|_{\Cop_{p,q} (u,v)} : = \big\| \|f\|_{p,v,(\cdot,b)} \big\|_{q,u,(a,b)} < \infty \bigg\},
		\end{align*}
		respectively.
	\end{defi}
	
	Note that if $0<\|u\|_{q,(t,b)}<\infty,~~ t>0,$ then $\Ces_{p,q} (u,v)$ is a quasi-normed vector space, also if $0<\|u\|_{q,(a,t)}<\infty,~~ t>0,$ then $\Cop_{p,q} (u,v)$ is a quasi-normed vector space.
	
	Many function spaces from the literature, in particular from harmonic analysis, are covered by the spaces $\Ces_{p,q}(u,v)$ and
	$\Cop_{p,q}(u,v)$. Let us only mention the Beurling algebras $A^p$
	and $A^*$, see \cite{gil1970, johnson1974, belliftrig}.
	
	The function spaces $C$ and $D$ defined by Grosse-Erdmann in \cite{Gr:98} are connected to our definition as follows:
	$$
	\Ces_{p,q}(u,v) = C(p,q,u)_v \qq \mbox{and} \qq \Cop_{p,q}(u,v) =
	D(p,q,u)_v.
	$$
	
	The classical Ces\`{a}ro function spaces $\Ces(p):=\Ces_{1, p}(x^{-1}, 1)$ were 
	introduced in 1970 by Shiue \cite{Sh:70} and were further investigated in
	\cite{Hass-Huss73} and \cite{SY-Zhang-Le:87}. Ces\`{a}ro and Copson function spaces, despite being defined correspondingly to Ces\`{a}ro sequence spaces, have not gotten as much attention as their sequence counterparts for a long time. In fact, there is extensive literature available on 
	different topics studied in Ces\`{a}ro sequence spaces as, for
	instance, in 
	\cite{Jag-74,Be:96, CuiPluc,cuihud1999,cuihud2001,cuihudli,cmp,chencuihudsims}. However, recently in a series of papers 
	\cite{astasmal2009,astasmal2010,astasmal2011,as,asmal13,asmal12},
	Astashkin and Maligranda started to study the structure of
	Ces\`{a}ro function spaces. In \cite{astasmal2009}, 
	they investigated dual spaces for $\Ces (p)$ for $1 < p < \infty$ (for further details on the classical Ces\`{a}ro spaces see \cite{asmalsurvey}).  
	
	The weighted  Ces\`{a}ro $C_{p,u}$ spaces studied by  Kami{\'n}ska and Kubiak in \cite{Ka-Ku12} are also a special case of our space  $ \Ces_{1,p}(u,{\bf 1}):=\Ces_{p,u} $.  In \cite{Ka-Ku12}, they computed the dual norm of the Ces\`{a}ro
	function space $\Ces_{p, u}$, generated by $1 < p < \infty$ and an
	arbitrary positive weight $u$.
	
	In \cite[Theorem 21.1]{Be:96} Bennett observed that the
	classical Ces\`{a}ro function space and the classical Copson
	function space $\Cop(p):=\Cop_{1,p}({\bf 1}, x^{-1})$ coincide for $p > 1$. He also derived estimates for the norms of the corresponding inclusion operators. The same result, with different estimates, is due to Boas \cite{Bo:70}, who, in
	fact, obtained the integral analog of the Askey-Boas Theorem \cite[Lemma 6.18]{Bo:67} and \cite{Hass-Huss73}. These
	results are generalized in \cite{Gr:98} using the blocking technique.

	If $X$ and $Y$ are two (quasi-) Banach spaces of measurable functions on $(a,b)$, we say that $X$ is \textit{embedded} into $Y$, and write $X \hookrightarrow Y$, if there exists a constant $c \in (0, \infty)$ such that $\|f\|_Y \leq c \|f\|_X$ for all $f \in X$. The \textit{optimal} (smallest possible) such $c$ is then $\|\Id\|_{X\rightarrow Y}$, where $\Id$ is the identity operator. A function $f$ is called a \textit{pointwise multiplier} from $X$ to $Y$ if the pointwise product $fg$ belongs to $Y$ for each $g\in X$. The space of all such multipliers $M(X,Y)$ becomes a quasi-normed space when endowed with the functional
	\begin{equation*}
		\|f\|_{M(X,Y)} = \sup_{g \nsim 0} \left\{ \frac{\|fg\|_Y}{\|g\|_X}\right\},
	\end{equation*}
	in which $g \nsim 0$ means that $g$ is not equimeasurable to the zero function (that is, $g$ is not a.e.~equal to zero).
	
	Pointwise multipliers between  $\Ces_{1, p}(x^{-1}, 1)$ and $\Cop_{1,q}(1, x^{-1})$ is given when $1 < q \leq p \leq \infty$ in \cite{kolema14} in 2019.
	
	If $f$ is a weight, $\|\cdot\|_{Y}\colon\mathcal M\to[0,\infty]$ is a functional and $Y\subset\mathcal M$ is given by
	\begin{equation*}
		Y=\{f\in\mathcal M:\|f\|_{Y}<\infty\},
	\end{equation*}
	then we can define the \textit{weighted space} $Y_f = \{g\in\mathcal M\colon fg \in Y\}$ and $\|g\|_{Y_f} = \|fg\|_Y$. Then one clearly has
	\begin{equation}\label{multiplier-embedding}
		\|f\|_{M(X,Y)} = \sup_{g \nsim 0} \left\{ \frac{\|g\|_{Y_f}}{\|g\|_X}\right\} = \|\Id\|_{X\rightarrow Y_f}.
	\end{equation}
	Therefore, the problem of finding pointwise multipliers between function spaces reduces to the characterization of the embeddings between weighted function spaces. 
	
	Our principal goal in this survey paper is to present the characterization of the embeddings between weighted Ces\`{a}ro and Copson function spaces, that is, the embeddings
	\begin{align}
		\Ces_{p_1,q_1}(u_1,v_1) & \hra \Ces_{p_2,q_2}(u_2,v_2), \label{mainemb1}\\
		\Cop_{p_1,q_1}(u_1,v_1) & \hra \Ces_{p_2,q_2}(u_2,v_2), 
		\label{mainemb2}\\
		\Cop_{p_1,q_1}(u_1,v_1) & \hra \Cop_{p_2,q_2}(u_2,v_2), \label{mainemb3}\\
		\Ces_{p_1,q_1}(u_1,v_1) & \hra \Cop_{p_2,q_2}(u_2,v_2).
		\label{mainemb4}
	\end{align}

	Note that by the change of variables $t \mapsto a+b-t$, in case when $b<\infty$ and $t \mapsto a+ \frac{1}{t-a}$, when $b=\infty$,  it is easy to see that \eqref{mainemb3} is
	equivalent to the embedding
	$$
	\Ces_{p_1,q_1}(\tilde{u}_1,\tilde{v}_1) \hra
	\Ces_{p_2,q_2}(\tilde{u}_2,\tilde{v}_2),
	$$
	where $\tilde{u}_i (t) := u_i\big(a+b-t\big)$,
	$\tilde{v}_i (t) :=  v_i\big(a+b-t\big)$,
	$i=1,2$, $t \in(a,b)$, when $b<\infty$, and $\tilde{u}_i (t) :=(t-a)^{ -2/q_2} u_i\big(a+\frac{1}{t-a}\big)$,
	$\tilde{v}_i (t) :=(t-a)^{-2/p_2} v_i\big(a+\frac{1}{t-a}\big)$, $i=1,2$, $t \in(a,\infty)$ when $b=\infty$. The same argument shows that \eqref{mainemb4} is
	equivalent to the embedding
	$$
	\Cop_{p_1,q_1}(\tilde{u}_1,\tilde{v}_1) \hra
	\Ces_{p_2,q_2}(\tilde{u}_2,\tilde{v}_2).
	$$
	This note allows us to concentrate our attention
	on the characterization of \eqref{mainemb1} and \eqref{mainemb2}.
	
	It is worth to mention that in view of \eqref{multiplier-embedding}, the characterizations of \eqref{mainemb1}
	- \eqref{mainemb2} makes it possible to characterize the pointwise multipliers between weighted Ces\`{a}ro and Copson function
	spaces which was previously noted as difficult to treat in \cite[p. 30]{Gr:98}.
	
	Note that the  embedding \eqref{mainemb2} was characterized in particular cases:
	
	{\rm (i)} $(a,b)=(0,\infty)$, 
	$p_1 = p_2 = 1$, $q_1 = q_2 = p > 1$, $u_1(t) = t^{\beta p - 1}$, $u_2(t) = t^{-\alpha p - 1}$, $v_1(t) = t^{-\beta - 1}$, $v_2(t) = t^{\alpha - 1}$, $t > 0$, where $\alpha > 0$ and $\beta > 0$, in \cite[p. 61]{Bo:70};
	
	{\rm (ii)} $(a,b)=(0,\infty)$,
	$p_1 = p_2 = 1$, $q_1 = p$, $q_2 = q$,
	$u_1(t) = v(t)$, $u_2(t) = t^{-q}w(t)$, $v_1(t) = t^{-1}$, $v_2 (t) = 1$, $t > 0$, where $0 < p \leq \infty$, $1 \leq q \leq \infty$ and $v,\,w$ are weight functions on $(0,\infty)$, in \cite[Theorem~2.3]{Ca-Go-Ma-Pi:08}, 
	
	{\rm (iii)} $(a,b)=(0,\infty)$,  $0<p_1,q_1<\infty$, $0 <p_2\leq q_2<\infty$, in \cite{GMU-CopCes, Go-Mu-Un:19},
	
	{\rm (iv)}  $(a,b)=(0,\infty)$, $p_1 = \infty$, $0<q_1<\infty$, $0 <p_2\leq q_2<\infty$,  in  \cite{Un:21}.
	
	The complete characterization of \eqref{mainemb2} without any restrictions on parameters or weights was given recently in \cite{GPU-JFA}.
	
	On the other hand, the characterization of the embedding \eqref{mainemb1} was given in the case when $(a,b) = (0,\infty)$, $0<p_1,q_1<\infty$, $0 <p_2\leq q_2<\infty$, in \cite{Un:20} by using duality in weighted Lebesgue spaces. By duality techniques, the embeddings \eqref{mainemb1}--\eqref{mainemb4} can be reduced to the weighted iterated inequalities (see, e.g., \cite{GKPS, GMU-CopCes, GMU-cLMLM, Un:20, Un:21}). 
	One of the earliest treatments of weighted inequalities involving iteration of operators was most likely carried out in \cite{Go-Mu-Pe:12},  \cite{Go-Mu-Pe:13} and \cite{Bu-Oi:13}. It is quite difficult to mention every contribution and provide a thorough history here, but interested readers can see
	\cite{MR-Saw89, BloKer91, Oinarov94, Step94, Lai99,  Oi-Ka:08, GogStep2013, Prok2016, Oi-Ka:15, Be-Or:17, Go-Mu:17, Ka:19, kal-oi:23,     KrePick-CC}. Recently, inequalities involving Hardy-Hardy or Hardy-Copson iteration were fully characterized in  \cite{Go-Mi-Pi-Tu-Un:21, Go-un:24} by eliminating duality approaches and hence removing numerous constraints that arise in previous works. 
	
	Furthermore,  in \cite[Theorem~2.1]{Go-Un:NS}, it was shown that the embeddings between local Morrey-type spaces are equivalent with some particular cases of \eqref{mainemb1}, however at that time the full characterization of the reduced inequalities were unknown. 
	
	In this paper, by avoiding duality techniques, we will present the characterization of \eqref{mainemb1} without any restrictions on the parameters or weight functions.
	
	\section{Embeddings between weighted  Ces\`{a}ro Spaces }\label{S:main-results}
	
	In this section we will focus on \eqref{mainemb1}. 
	
	\begin{rem}\label{Rem-pp} \cite[Lemmas 3.4-3.5]{Go-Mu:17}
		Note that $\Ces_{p,p}(u,v)$ and $\Cop_{p,p}(u,v)$ coincide with some weighted Lebesgue spaces. 
		Indeed, 
		\begin{itemize}
			\item if $w(x) =v(x)  \|u\|_{p,(x,b)}$, then $\Ces_{p,p}(u,v) = L_p(w)$.

			\item if $w(x) = v(x) \|u\|_{p,(a,x)}$, then $\Cop_{p,p}(u,v) = L_p(w)$.
		\end{itemize} 
	\end{rem}
	In view of Remark~\ref{Rem-pp} we have the following:
	\begin{cor} Let $0< p_1, p_2, q_1, q_2 \leq \infty$ and $u_1, u_2, v_1, v_2 \in \mathcal{W}(a,b)$.
		\begin{itemize}
			\item[(i)] The embedding $$\Ces_{p_1,p_1}(u_1,v_1)  \hra \Ces_{p_2,q_2}(u_2,v_2)$$ 
			is reduced to the embedding
			\begin{equation} \label{cor1}
				L_{p_1}(w_1) \hra \Ces_{p_2,q_2}(u_2,v_2)
			\end{equation}
			with $w_1(x)= v_1(x)  \|u_1\|_{p_1,(x,b)}$.
			\item[(ii)] The embedding $$\Ces_{p_1,q_1}(u_1,v_1)  \hra \Ces_{q_2,q_2}(u_2,v_2)$$ 
			is reduced to the embedding
			\begin{equation}\label{cor2}
				\Ces_{p_1,q_1}(u_1,v_1)\hra L_{q_2}(w_2)
			\end{equation}
			with $w_2(x)= v_2(x)  \|u_2\|_{q_2,(x,b)}$.
			
			\item[(iii)] The embedding $$\Cop_{p_1,p_1}(u_1,v_1)  \hra \Ces_{p_2,q_2}(u_2,v_2)$$
			is reduced to the embedding
			\begin{equation}\label{cor3}
				L_{p_1}(w_1)  \hra \Ces_{p_2,q_2}(u_2,v_2)
			\end{equation}
			with 
			$w_1(x) = v_1(x) \|u\|_{p_1,(a,x)}.$
			
			\item[(iv)]The embedding $$\Cop_{p_1,q_1}(u_1,v_1)  \hra \Ces_{q_2,q_2}(u_2,v_2)$$ 
			is reduced to the embedding
			\begin{equation}\label{cor4}
				\Cop_{p_1,q_1}(u_1,v_1)  \hra L_{q_2}(w_2)
			\end{equation} 
			with 
			$w_2(x) = v_2(x) \|u\|_{q_2,(a,x)}$.
		\end{itemize}
	\end{cor}
	
	The characterization of the inequalities  \eqref{cor1}--\eqref{cor4} are given in \cite{GMU-CopCes}. These inequalities are equivalent to the well-established weighted Hardy and weighted reverse  Hardy inequalities. It would be impossible to mention every contribution here. The development and characterization of weighted Hardy inequalities can be seen in the following books \cite{op-ku90} and \cite{ku-pe-sa17}; also, the history of Hardy inequalities can be found in \cite{ku-pe-ma07}. 
	
	The reverse Hardy inequalities on the other hand, are studied by Beesack~\cite{Bee:61}, Beesack and Heinig~\cite{beehei81},  Prokhorov~\cite{prok04}, Evans,  Gogatishvili and Opic~\cite{evgoop:08},   Mustafayev and \"{U}nver~\cite{Musunv15}.
	
	\begin{rem} The characterization of the embeddings \eqref{cor2} and \eqref{cor4} allow us to give a characterization of the associate spaces of weighted Ces\`{a}ro and Copson function spaces (see Theorems~3.11-3.12 from \cite{GMU-CopCes}).
	\end{rem}
	
	Our aim is to characterize the weights $u_1, u_2, v_1, v_2$ for which the inequality
	\begin{align}\label{Ces-Ces}
		&\bigg(\int_a^b \bigg(\int_a^t f(s)^{p_2} v_2(s)^{p_2} ds \bigg)^{\frac{q_2}{p_2}} u_2(t)^{q_2} dt \bigg)^{\frac{1}{q_2}} \nonumber\\
		& \hspace{3cm} \leq c  \bigg(\int_a^b \bigg(\int_a^t f(s)^{p_1} v_1(s)^{p_1} ds \bigg)^{\frac{q_1}{p_1}} u_1(t)^{q_1} dt \bigg)^{\frac{1}{q_1}}
	\end{align}
	holds for every measurable non-negative function $f$ on $(a,b)$. 
	
	To this end, we will consider a more manageable version of~\eqref{Ces-Ces}. It is clear that when we swap out $f$ for $(fv_1)^{p_1}$, the set of competing functions in~\eqref{Ces-Ces} remains the same. Raising~\eqref{Ces-Ces} to $p_1$, defining new parameters $p,q,r$ by
	\begin{equation*}
		r=\frac{p_2}{p_1}, \,\, q= \frac{q_2}{p_1}, \,\, \text{and} \,\, p= \frac{q_1}{p_1}
	\end{equation*}
	and introducing new weights $u,v,w$ by
	\begin{equation*}
		u=u_2^{q_2}, \,\, v= v_1^{-p_2} v_2^{p_2}, \,\, \text{and} \,\, w= u_1^{q_1},
	\end{equation*}
	We obtain an equivalent form of~\eqref{Ces-Ces}, that is,
	\begin{equation}\label{ces-ces-main}
		\bigg(\int_a^b \bigg(\int_a^t f(s)^r v(s) ds \bigg)^{\frac{q}{r}} u(t) dt \bigg)^{\frac{1}{q}} \leq C \bigg(\int_a^b \bigg(\int_a^t f(s) ds \bigg)^{p} w(t) dt\bigg)^{\frac{1}{p}}.
	\end{equation}
	If $c$ and $C$ are the best constants respectively in~\eqref{Ces-Ces} and~\eqref{ces-ces-main}, then $c^{p_1} = C$.
	
	Therefore, we can focus on inequality~\eqref{ces-ces-main}. We aim to give a characterization without any restrictions on parameters and weights. To accomplish this goal, we will avoid duality techniques and construct a new discretization and anti-discretization approach. 
	
	\begin{rem} It was shown in \cite[Lemma~1]{Un:20} that \eqref{ces-ces-main} holds only for trivial functions if $r>1$. Hence, the assumption that $0 < r \leq 1$ in Theorem~\ref{T:ces ces main} is not a restriction.
	\end{rem}
	
	For $0\leq a < b \leq \infty$, we define
	\begin{align}\label{Vr}
		V_r(a, b) =
		\begin{cases}
			\bigg(\int_a^b v^{\frac{1}{1-r}} \bigg)^{\frac{1-r}{r}} & \text{if} \quad  0<r<1,
			\\
			\esup\limits_{s \in (a, b)} v(s) & \text{if} \quad  r=1.
		\end{cases}
	\end{align}
	
	Our main result is the following. 
	
	\begin{thm}\label{T:ces ces main}
		Let $0 < r \leq 1$, $0 < p, q < \infty$ and assume that $u,v,w$ are weights on $(a,b)$ such that $0 < \int_x^b w<\infty$ for all $x\in (a,b)$. Then \eqref{ces-ces-main} holds for all $f\in \mathcal{M}^+(a,b)$ if and only if
		
		\rm{(i)} either  $p \leq r \leq 1 \leq q$ and $C_1< \infty$, where
		\begin{equation*}
			C_1 :=  \esup_{x\in (a,b)} \bigg(\int_x^b w\bigg)^{-\frac{1}{p}} \esup_{t\in (x,b)} \bigg( \int_t^b u \bigg)^{\frac{1}{q}}  V_r(x,t),
		\end{equation*}
		moreover, the best constant $C$ in \eqref{ces-ces-main} satisfies $C \approx C_1$,
		
		\rm{(ii)} or  $p\leq q <1$, $p \leq r\leq 1$ and $C_2 < \infty$, where
		\begin{equation*}
			C_2 :=  \sup_{x\in (a,b)} \bigg(\int_x^b w\bigg)^{-\frac{1}{p}} \bigg( \int_x^b \bigg( \int_t^b u \bigg)^{\frac{q}{1-q}} u(t) V_r(x,t)^{\frac{q}{1-q}} dt \bigg)^{\frac{1-q}{q}},
		\end{equation*}
		moreover, the best constant $C$ in \eqref{ces-ces-main} satisfies $C \approx C_2$,
		
		\rm{(iii)} or  $r < p \leq q$, $r \leq 1 \leq q$, and $\max\{C_1, C_3\} < \infty$, where
		\begin{equation*}
			C_3 :=  \sup_{x \in (a,b)} \bigg(\int_{x}^b u\bigg)^{\frac{1}{q}} \bigg(\int_a^{x} \bigg(\int_t^b w \bigg)^{-\frac{p}{p-r}} w(t) V_r(t,x)^{\frac{pr}{p-r}} dt \bigg)^{\frac{p-r}{pr}},
		\end{equation*}
		moreover, the best constant $C$ in \eqref{ces-ces-main} satisfies $C \approx C_1+ C_3$,
		
		\rm{(iv)}  or $r < p \leq q < 1$, 
		and $\max\{C_2, C_3\} < \infty$, moreover, the best constant $C$ in \eqref{ces-ces-main} satisfies $C \approx C_2 + C_3$,
		
		\rm{(v)} or $q < p \leq r \leq 1$, and $\max\{C_4, C_5\} < \infty$, where
		\begin{equation*}
			C_4 := \bigg(\int_a^{b} \bigg(\int_x^b u\bigg)^{\frac{q}{p-q}} u(x)  \esup_{t\in (a, x)}  \bigg(\int_t^b w \bigg)^{-\frac{q}{p-q}} V_r(t,x)^{\frac{pq}{p-q}} dx \bigg)^{\frac{p-q}{pq}},
		\end{equation*}
		and
		\begin{align*}
			C_5 & :=  \bigg(\int_a^b  w(x) \sup_{y\in (a,x)} \bigg(\int_y^b w\bigg)^{-\frac{p}{p-q}} \\
			& \hspace{2cm} \times \bigg( \int_{y}^x \bigg( \int_t^x u \bigg)^{\frac{q}{1-q}} u(t) V_r(y,t)^{\frac{q}{1-q}} dt \bigg)^{\frac{p(1-q)}{p-q}} dx \bigg)^{\frac{p-q}{pq}},
		\end{align*}
		moreover, the best constant $C$ in \eqref{ces-ces-main} satisfies $C \approx C_4 + C_5$,
		
		\rm{(vi)} or  $q < p$, $q < 1$, $r < p$, $r \leq 1$, and $\max\{C_1, C_5, C_6\} < \infty$, where
		\begin{align*}
			C_6 &:= \Bigg(\int_a^{b} w(x) \esup_{ y\in(a, x)} \bigg(\int_y^{b} w\bigg)^{-1} \Bigg(\int_y^{x} 
			\bigg(\int_t^b u\bigg)^{\frac{q}{p-q}} u(t)dt\\
			& \hskip+3cm \times \bigg(\int_a^y  \bigg(\int_t^b w \bigg)^{-\frac{p}{p-r}} w(t)V_r(t,y)^{\frac{pr}{p-r}}dt\bigg)^{\frac{(p-r)q}{(p-q)r}} dx \Bigg)^{\frac{p-q}{pq}},
		\end{align*}
		moreover, the best constant $C$ in \eqref{ces-ces-main} satisfies $C \approx C_1 + C_5 + C_6$.
		
		\rm{(vii)} or   $r \leq 1 \leq q < p$, and $\max\{C_1, C_6, C_7\} < \infty$, where
		\begin{align*}
			C_7 &= \Bigg(\int_a^{b}  w(x) \esup_{y\in(a,x) } \bigg(\int_y^b w\bigg)^{-\frac{p}{p-q}}  \bigg(\esup_{t\in(y,x) }  \bigg( \int_t^{x} u \bigg)^{\frac{p}{p-q}} V_r(y,t) ^{\frac{pq}{p-q}}\bigg) dx \Bigg)^{\frac{p-q}{pq}},
		\end{align*}
		moreover, the best constant $C$ in \eqref{ces-ces-main} satisfies $C \approx C_1 + C_6 + C_7$.
	\end{thm}


	%
	
\end{document}